\documentclass[a4paper,12pt]{article}
\usepackage[cp1251]{inputenc}
\usepackage{amssymb,amsmath,amsthm}
\usepackage[russian]{babel}
\usepackage{amssymb,amsmath,amsthm}
\usepackage{indentfirst}
\usepackage{hyperref,hypbmsec}
\usepackage[pdftex]{graphicx}
\numberwithin{equation}{section}
\newtheorem{Th}{\hskip\parindent Theorem}[section]
\newtheorem{De}{\hskip\parindent Definition}[section]
\newtheorem{Le}{\hskip\parindent Lemma}[section]

\newtheorem{Zam}{\hskip\parindent Remark}[section]

\newcounter{propet}

\renewcommand{\le}{\leqslant}\renewcommand{\ge}{\geqslant}

\begin{document}
 \selectlanguage{russian}

\begin{Large}
\centerline{\bf{Methods for estimating of continuants, corrected version.}}
\centerline{ I. D. Kan}
\end{Large}
\vskip+1.0cm

{\bf Abstract.} 
The purpose of \cite{0} was as follows.
``We consider special sets of continuants  which occur in applications.
For these sets  we solve the problem of finding maximal and minimal continuants.
There are several  methods for finding extremum
such as  thee: mathod of basic substitutions (section 4.1),
method of quadratic irrationalities  (section 4.2), 
majorizing inequalities (section 4.3)
and 
unit variation (section 4.4).''
Unfortunately, the text of \cite{0} was littered with numerous typos, making it difficult to read. 
This defect is corrected in the present work. 
In addition, some of the statements are now provided with more detailed evidence.
\section{  Introduction.}

\tableofcontents

For integer 
$t\geq0$ and positive integers 
$a_1,a_2,\ldots,a_t$ the continuant
$\langle a_1,a_2,\ldots,a_t\rangle $
is defined inductively by equalities $\langle 
\varnothing\rangle=1$, $\langle a_1\rangle=a_1$,
\begin{equation} 
\label{1} 
\langle a_1,a_2,\ldots,a_{j+1}\rangle
=
a_{j+1}\langle 
a_1,a_2,\ldots,a_j\rangle
+
\langle a_1,a_2,\ldots,a_{j-1}\rangle
,
\end{equation} 
$j=1,2,\ldots,t-1$. 
For a continued fraction
\begin{equation}
\label{cont.fraction}
[a_0;a_1,a_2,\ldots,a_t]
=
a_0 +
\cfrac{1}{a_1+\cfrac{1}
{a_2+{\atop\ddots\,
\displaystyle{
+
\cfrac{1}{a_t}}}}}
,
\end{equation}
we know (see, \cite{1}) that
\begin{equation} 
\label{2} 
[a_0;a_1,a_2,\ldots,a_t]
=
\frac{\langle a_0,a_1,\ldots,a_t\rangle}
{\langle a_1,a_2,\ldots,a_t\rangle}
.
\end{equation} 
We suppose the continued fraction of the empty set to be equal to zero.

Usually it is necessary to consider kontinualts with elements boulded from above by 
$n$ (see \cite{7}-\cite{9-10}). 
In certain problems in number theory  (see 
\cite{8}) the following condition occurs:
$$
a_1+a_2+\ldots+a_t=S, \,\,\,\, S/t\ll{const}.
$$
Such restrictions define certain sets of continuants.
The problem arises is to find the maximal and the minimal value of continuant
over these sets.
Some results related to this problem were obtained in 
\cite{2,11,12}.
In the present paper we  obtain further reuslts.

The author is grateful to Prof. N. Moshchevitin for setting the problem and  for  idea of   method of section 4.4.

\section{  Notation and parameters.}\,\,

For $a,b\in\mathbb{N}\cup\{0\}$ put

$$
a^{\{b\}}
=
\begin{cases}
\underbrace{\{a,a,\ldots,a\}}\limits_{\mbox{$b$ \text{elements}}},& \mbox{if $b\geq1$},\cr \varnothing,  &\mbox{ if $b=0$}
\end{cases}
.
$$

For a collection of finite sequences $A_1,A_2,\ldots
, A_k$ 
by $\{A_1,A_2,\ldots A_k\}$
we define the sequence
of all  consecutive elements
of $A_1,A_2,\ldots A_k$

Consider
$f$ natual numbers $h_1<h_2<\ldots<h_f$
with multiplicities
$p_1,p_2,\ldots,p_f\geq1$. 
Let's put it 
$$
W_f=W_f(\bar{h},\bar{p})=\left\{(a_1,a_2,\ldots,a_t):a_i\in\{h_1,h_2,\ldots,h_f\},\,\,|\{i:a_i=h_j\}|=p_j\right\},
$$
where $t=p_1+\ldots+p_f$. 
Consider a partition of each
 $p_1,\ldots,p_f$ into a sum
\begin{equation} 
\label{3} 
p_j=l_j+r_j 
\end{equation} 
of two nonnegative $l_j$, $r_j$ and put
$$
D_f=D_f(\bar{l},\bar{r})=\left\{
h_f^{\{l_f\}},h_{f-1}^{\{l_{f-1}\}},\ldots,h_2^{\{l_2\}},
h_1^{\{l_1\}},h_1^{\{r_1\}},h_2^{\{r_2\}},\ldots,h_f^{\{r_f\}}
\right\}.
$$
Define
$$
V_f=V_f(\bar{h},\bar{p})=\left\{(a_1,a_2,\ldots,a_t)\in
W_f(\bar{h},\bar{p}):a_1=h_1\right\}.
$$
Moreover for $S,t,n\in\mathbb{N}$ put
$$
U(S,t,n)=\left\{(a_1,a_2,\ldots,a_t):a_i\in\{1,2,\ldots,n\},a_1+a_2+\ldots+a_t=S\right\},
$$
$$
U_n(S)=\bigcup\limits_tU(S,t,n),\,\,\,
U(S,t)=\bigcup\limits_nU(S,t,n)
.
$$
Here by $[\xi],\{\xi\} $
we denote integer part and fractional part for 
$\xi\in\mathbb{R}$, correspondingly.

From the set  $W_f(\bar{h},\bar{p})$ we extract two monotone sequences
\begin{equation} 
\label{4} 
b_0\leq{b_1}\leq\ldots\leq{b_{t-1}}    
\end{equation} 
and
\begin{equation} 
\label{4a} 
c_1\geq{c_2}\geq\ldots\geq{c_t} 
.
\end{equation}

For $j=0,1,\ldots , t$ and $\nu=1,2,\ldots,2[t/2]$ put
$$n_j=\left\{
\begin{array}{cl}
b_j,\,\,\mbox{if}\,\, j\equiv0\pmod2,\\
c_j,\,\,\mbox{if} \,\, j\equiv1\pmod2,
\end{array}
\right. m_\nu=\left\{
\begin{array}{cl}
c_\nu,\,\,\mbox{if} \,\,\nu\equiv0\pmod2,\\
b_\nu,\,\,\mbox{if} \,\, \nu\equiv1\pmod2.
\end{array}
\right.
$$

For $m\in\mathbb{N}\cup\{0\}$ put $m_-=[m/2]$, $m_+=m-m_-$. 
For
$x\in\mathbb{Z}$,$n,z\in\mathbb{N}$ define
\begin{equation} 
\label{5} 
N_z(x)=\left\{
\begin{array}{lc}
\left\{n^{\{x\}} , z\right\},&\mbox{ if x}\geq0,\\
\left\{1^{\{-1-x\}},z,1\right\},&\mbox{ if x}<0,
\end{array}
\right.
\end{equation}
\begin{equation} 
\label{6} 
T_z(m,x)=
\langle\underbrace{1,n,1,n,\ldots,1,n}_{2 m_+ \mbox
{elements}} ,N_z(x),\underbrace{n,1,n,1\ldots,n,1}_{2 m_-
\mbox{elements}}\rangle  
.
\end{equation}

For $S,n\in\mathbb{N}$ define
\begin{equation} 
\label{7} 
S_0^{(n)}
=
(n+1)\left\{
\frac{S-1}{n+1}
\right\}+1
,\,\,\,\,\,\,
S_1^{(n)}
=
n\left\{\frac{S}{n}\right\}
. 
\end{equation}
We define a subset $P(S)$ of the set $\{1,2,\ldots,n-1\}$
by the condition 
\begin{equation} 
\label{8} 
P(S)=\left\{
\begin{array}{ll}
\{1,2,\ldots,n-1\}, & \mbox{if } S\geq{n^2-1},\\
\{S_0^{(n)},S_0^{(n)},\ldots,S_1^{(n)}\}, & \mbox{if } S<n^2-1\, \mbox{and } S_0^{(n)}\leq S_1^{(n)},\\
\{1,2,\ldots,S_1^{(n)}\} \bigcup\,
\{S_0^{(n)},S_0^{(n)}+1,\ldots,n-1\}  & \mbox{in the opposite case.}
\end{array}
\right. 
\end{equation}

\section{  Main results.}\,\,

\begin{Th}\label{t1} \cite{2}
The maximal continuant over the set 
$V_f(\bar{h},\bar{p})$ satisfies 
\begin{equation} 
\label{9} 
\max\limits_{V_f}\langle a_1,a_2,\ldots,a_t\rangle
=
\langle b_0,b_1,\ldots,b_{t-1}\rangle
,
\end{equation}
where the maximum is taken over all permutations of the sequence $ \{a_1,a_2,\ldots,a_t\}$.
\end{Th}

\begin{Th}\label{t2}The maximal continuant over the set
$W_f(\bar{h},\bar{p})$ satisfies
\begin{equation} 
\label{10} 
\max_{W_f}\langle a_1,a_2,\ldots,a_t\rangle
=
\langle D_f(\bar{l},\bar r)\rangle
, 
\end{equation}
here integers $l_j$ and $r_j$ $(j=1,2,...,f)$ from (\ref{3}) are such that
\begin{equation} 
\label{11} 
1=\left\{
\begin{array}{cll}
l_j, &\mbox {if} \,\,j\equiv f &\pmod {2},\\
r_j, &\mbox {if} \,\,j\equiv f-1& \pmod {2}
.
\end{array}
\right. 
\end{equation}
\end{Th}

\begin{Th}\label{t3}
The minimal continuant over the set
$W_f(\bar{h},\bar{p})$ satisfies
\begin{equation} 
\label{12} 
\min_{W_f}\langle a_1,a_2,\ldots,a_t\rangle
=
\langle n_0,n_1,\ldots,n_\nu,
m_\mu,m_{\mu-1},\ldots,m_2,m_1\rangle
, 
\end{equation}
where $\,\,\,\,\nu
=
t_+-1$, $\mu=t_-$.
\end{Th}

\begin{Th}\label{t4}
The maximal continuant over the set
  $U_n(S)$  satisfies
\begin{equation} 
\label{13} 
\max_{U_n(S)}
\langle 
a_1,a_2,\ldots,a_t
\rangle 
=
\langle \underbrace{1,1,\ldots,1}_{S\text{
elements}}\rangle 
.
\end{equation}
\end{Th}

\begin{Th}\label{t5}
For $2\leq t\leq S$  the maximal continuant over the set $U(S,t)$ 
satisfies
\begin{equation} 
\label{14} 
\max_{U(S,t)}
\langle 
a_1,a_2,\ldots,a_t
\rangle 
=
\left\{
\begin{array}{lc}
\langle h_1^{\{ t\}}\rangle
, & \mbox{if}\,\, 
S\equiv 0 \pmod{t},\\
\langle h_2,h_1^{\{ c\}}
,
h_2^{\{ d-1\}}\rangle
,& 
\mbox{ --- in the opposite case,}
\end{array}
\right. 
\end{equation}
where
\begin{equation} 
\label{15} 
h_1=[S/t], \, \,\,
h_2=h_1+1, \, c=t\{-S/t\}, \, d=t\{S/t\}
. 
\end{equation}
\end{Th}

\begin{Th}\label{t6}
For $2\leq t\leq S\leq nt$ the minimal continuant over the set $U(S,n,t)$ satisfiies
\begin{equation} 
\label{16} 
\min_{U(S,n,t)}\langle a_1,a_2,\ldots,a_t\rangle
=
T_z(m,x), 
\end{equation}
where
\begin{equation} 
\label{17} 
z\in\{1,2,\ldots,n-1\}, \,\,\,\, z\equiv S-t+1\pmod{(n-1)}
,  
\end{equation}
\begin{equation} 
\label{18} 
x
=
\frac{2(S-t+1-z)
}{n-1} 
+1-t
,
\,\,\,\,\,
m
=
(t-|x|-1)/2
.
\end{equation}
\end{Th}

\begin{Th}\label{t7}
For $n\geq2$, $S\geq2n+2$ the minimal continuant over  the set $U_n(S)$ satisfies
\begin{equation} 
\label{19} 
\min_{U_n(S)}\langle a_1,a_2,\ldots,a_t\rangle
=
\min_{z\in P^{*}(S)}T_z(m,x)
, 
\end{equation}
where
\begin{equation} 
\label{20} 
0\leq x\leq n,  \, x\equiv z-S\pmod{(n+1)}, \,\,\,\, 
m=m(z)=\frac{S-z-nx}{n+1}
, 
\end{equation}
$$
P^{*}(S)=\{z\in P(S):\, m(z)\geq1\}.
$$
\end{Th}

{\bf Corolary 1.}\,\,{\it For $2\leq n\leq S-2$ one has
\begin{equation} 
\label{21} 
\min_{U_n(S)}\langle a_1,a_2,\ldots,a_t\rangle
\geq
\frac{\left(\sqrt[n+1]{\mu_n}\right)^{S-n+1}
}{e^2}
, 
\end{equation}
where 
\begin{equation} 
\label{22} 
\mu_n=\frac{n+2+\sqrt{n^2+4n}}{2}
.
\end{equation}
}

\begin{Zam}\label{z1} 
In  \cite{7}
for the minimum form Corollary 1 for $ n=4$ the authors use the lower bound of the form
$\left(\sqrt{2+10^{-6}}\right)^S$. 
The inequality  (\ref{21}) gives a lower bound of the form
$\left(\sqrt[5]{3+\sqrt{8}}\right)^S$. 
For large values of $S$
(for $S\geq350$) this bound is more precise as
$$
\sqrt[5]{3+\sqrt{8}}=1,424\ldots>1,414\ldots
=
\sqrt{2+10^{-6}}.
$$
\end{Zam}

\begin{Zam}\label{zz2} 
The \ref{t4} theorem, given here for 
"completeness'', should be considered part of mathematical folklore: it is quite obvious and is used in many works. 
For example, N. M. Korobov in 1997 in the work
\cite{kor97} formulated some algorithm related to continued fractions. 
The proof of this result (which he left to his readers) seems to rely in part on the \ref{t4} theorem.
\end{Zam}

\begin{Zam}\label{zz3} 
Statements of problems on calculating the maximum and minimum from the \ref{t2} and \ref{t3} theorems are available in the book by D. E. Knuth \cite{knuth1}. 
The solution of the second one coincides with the formulation of the \ref{t3} theorem.
\end{Zam}

\section{  Auxilary statements.}
\subsection{   Basic substitution method.}
 
 Let us consider the inequality
\begin{equation} 
\label{23} 
\langle u_1,u_2,\ldots,u_\alpha,v_1,v_2,\ldots,v_\beta,w_1,w_2,\ldots,w_\gamma\rangle
\leq
\langle u_1,u_2,\ldots,u_\alpha,v_\beta,v_{\beta-1},\ldots,v_2,v_1,w_1,w_2,\ldots,w_\gamma\rangle
.
\end{equation}

\begin{Le}\label{1A}
\cite{11}. 
For $v_1\neq v_\beta$, $u_\alpha\neq w_1$
the inequality (\ref{23}) 
is valid if and only if 
\begin{equation} 
\label{24} 
(v_1-v_2)(u_\alpha-w_1)<0
. 
\end{equation}
\end{Le}

Put
\begin{equation} 
\label{25} 
[\stackrel{\leftarrow}{U}]=[0;u_\alpha,u_{\alpha-1},\ldots,u_2,u_1], 
\end{equation}
\begin{equation} 
\label{26} 
[\overline{W}]=[0; w_1,w_2,\ldots,w_\gamma], 
\end{equation}
$$
[\overline{V}]=[0; v_1,v_2,\ldots,v_\beta], 
$$
$$
[\stackrel{\leftarrow}{V}]
=
[0;v_\beta,u_{\beta-1},\ldots,u_2,u_1],
$$
\begin{equation} 
\label{26a} 
{\bf a'}(U,V,W)
=
([\stackrel{\leftarrow}{U}]-[\overline{W}])
(v_1-v_\beta)
,
\end{equation}
\begin{equation} 
\label{26b} 
{\bf a}(U,V,W)
=
([\stackrel{\leftarrow}{U}]-[\overline{W}])
([\stackrel{\leftarrow}{V}]-[\overline{V}])
.
\end{equation}

In order to use Lemma \ref{1A} with
 $[\stackrel{\leftarrow}{U}] =0$ or $[\overline{W}]=0$
it is convenient to put  $u_\alpha=+\infty$ or
$w_1=+\infty$.
We suppose that adding the sigh ``infinity'' 
 does not change the value of the continuant or continued fraction's value, by the definition.
The sign ``infinity'' may appear in the inequality (\ref{24}).

\begin{Le}\label{1B} \cite{2}
 For   $v_1\neq v_\beta$ the inequality  (\ref{23})
is valid if and only if
one of two foollowing groups of conditions are valid:

the first group
$$
v_1>v_\beta, \,\, 
[\stackrel{\leftarrow}{U}]\geq[\overline{W}]
,
$$
or the second group
\begin{equation} 
\label{27} 
v_1<v_\beta, 
\end{equation}
\begin{equation} 
\label{28} 
[\stackrel{\leftarrow}{U}]\leq[\overline{W}]
. 
\end{equation}
In other words, the inequality (\ref{23}) holds if and only if
$$
{\bf a}(U,V,W)
\geq
0
;
$$
or, if and only if
$$
{\bf a'}(U,V,W)
\geq
0
.
$$
Moreover
the inequality (\ref{23})
turns into an equality if and only if
$[\stackrel{\leftarrow}{U}]
=
[\overline{W}]$, where
$[\stackrel{\leftarrow}{U}],[\overline{W}],$ 
${\bf a}(U,V,W)$, ${\bf a'}(U,V,W)$ are defined in (\ref{25}) --- (\ref{26b}).
\end{Le}

\begin{De}\label{d1} 
Suppose $\alpha+\beta+\gamma=t$,
\begin{equation} 
\label{29} 
\{a_1,a_2,\ldots,a_t\}=\{u_1,u_2,\ldots,u_\alpha,v_1,v_2,\ldots,v_\beta,w_1,w_2,\ldots,w_\gamma\}.
\end{equation}
Suppose 
$\Pi$ 
to be a substitution of indices
$\{1,2,\ldots,t\}$, which changes the sequence (\ref{29}) into the sequence
\begin{equation} 
\label{30} 
\left\{a_{\Pi(1)},a_{\Pi(2)},\ldots,a_{\Pi(t)}\right\}=\{u_1,u_2,\ldots,u_\alpha,v_{\beta},v_{\beta-1},\ldots,
v_2,v_1,w_1,w_2,\ldots,w_\gamma\}. 
\end{equation}
The we define 
$\Pi$
to be a basic substitution with the middle part
$(v_1,\ldots,v_\beta)$ or briefly $(\alpha+1,\alpha+\beta)$-substitution.
\end{De}

The simplest property is the symmetry of continuants (see \cite{11}):
\begin{equation} 
\label{31} 
\langle a_1,a_2,\ldots,a_t\rangle 
=
\langle a_t,a_{t-1},\ldots,a_2,a_1\rangle  
.
\end{equation}
Another simple property is the property of ``unit extraction'':
\begin{equation} 
\label{32} 
\langle a_1+1,a_2,\ldots,a_t\rangle 
=
\langle 1,a_1,a_2,\ldots,a_t\rangle 
. 
\end{equation}

\begin{De}\label{d2} 
A basic substitution
 $\Pi$
we define to be trivial with respect to the sequence
(\ref{29}),
if the continuant (\ref{29}) is equal to the continuant (\ref{30}), and this equality may be proved by means
(\ref{31}) and (\ref{32}) only.
\end{De}

{\bf Example 1.}\,\,
The substitution
 $\Pi:(2,4,5,1,1)\mapsto(2,5,4,1,1)$
is trivial as from (\ref{31}) and (\ref{32}) it follows that
$$
\langle 2,4,5,1,1\rangle =\langle 1,1,4,5,1,1\rangle =
\langle 1,1,5,4,1,1\rangle 
=
\langle 2,5,4,1,1\rangle 
.
$$

\begin{Zam}\label{z2} 
If a basic substitution does not change the value of a continuant then this substitution is trivial.
\end{Zam}

\begin{De}\label{d3} 
Consider a basic substitution $\Pi$ with the middle
$(v_1,v_2,\ldots,v_\beta)$.
Supposse the inequality (\ref{23}) to be strict one (not necessary strict one).
Then we define the basic substitution to be  increasing (non-decreasing) with respect to  the sequence 
(\ref{29}).
Analogously in the case of strict (not necessary strict)
 inequality opposite to 
(\ref{23}) we define 
$\Pi$ to be decreasing (non-increasing) with respect to the sequence
(\ref{29}).
\end{De}

Suppose that
$1\leq k<l\leq t$.
Then the sequence
$\{a_k,a_{k+1},\ldots,a_l\}$ we define to be a  $(k,l)$-fragment
of the sequence
 $\{a_1,a_2,\ldots,a_t\}$ with boundaries $a_k$ and
$a_l$. For $k\neq1$ or $l\neq t$ the boundaries $a_k$ and $a_l$
(correspondingly) we define as the proper boundaries.

\begin{Zam}\label{z3} 
If a basic substitution is increasing (decreasing) for a fragment of a sequence with a middle which does
not contain proper boundaries  of the fragment the this substitution is
increasing (decreasing) for the whole sequence itself.   
 
(In general, this property is no longer valid for trivial substitutions.)
\end{Zam}

\begin{Le}\label{l2}
For integers $r,h,g$ under the conditions  $r\geq 0, \,
0<h<g\leq+\infty$ one has
\begin{equation} 
\label{33} 
[0;h^{\{ r\}},g]\leq1/h. 
\end{equation}
Moreover  (\ref{33}) 
turns into an equality only in the case
 $r=1$, $g=+\infty$.
\end{Le}

\begin{Le}\label{l3}
The maximal continuant over the set
$W_f(\bar{h},\bar{p})$ 
is attained on a sequence 
$D_f(\bar{l},\bar{r})$ with elements satisfying  (\ref{3}).
\end{Le}

Proof. 

Suppose that the maximal continuant over the set
$W_f(\bar{h},\bar{p})$
is attained on the sequence
$D=\{a_1,\ldots,a_t\}$. Suppose that for a certain value of 
$\nu \in\{1,2,\ldots,t\}$ one has $a_\nu=h_1$. Consider
two basic substitutions
$\Pi$ and $\Phi$ acting on $(1,\nu)$- and $(\nu,t)$-fragments of $D$, correspondingly.
Suppose that these substitutions satisfy  $a_{\Pi(\nu)}=a_{\Phi(\nu)}=h_1$.
Then  $\Pi$ and $\Phi$  are non-increasing for the corresponding fragments.
Then from Theorem (\ref{t1} (Theorem (\ref{t1} is proved without application of Lemma (\ref{l3})
we know the srtucture of the fragments under consideration.
We take into account Remark (\ref{z3}). 
Lemma follows.

\begin{Le}\label{l4}
Suppose that elements of the sequence  $D_f(\bar{l},\bar{r})$
satisfy
\begin{equation} 
\label{34} 
l_f=1, \, r_f=0, \, r_{f-1}>0 
\end{equation}
Suppose that for this sequence the exist no 
increasing basic substitutions.
Then elements of the sequence  $D_f$ satisfy (\ref{11}).
\end{Le}

Proof.

Suppose that
 $f\geq3$ (otherwise there is nothing to prove).
Now we process induction in   $f$. 
By the assumption, the equality
(\ref{11}) is true for 
$j=f$. 
We will show that 
$r_{f-1}=1$ and $l_{f-2}>0$. 
Assume that it is not so.
Supoose that
 $r_{f-1}>1$ or $(r_{f-1}=1)\;\&(l_{f-2}=0)$. 
For both of these cases we can consider an 
$(l_{f-1}+2,t-1)$-basic substitution. 
Now we allpy Lemma \ref{1B}.
inequality (\ref{27}) is valid, by the construction.
We have also a strict inequality in 
(\ref{28}), by Lemma \ref{l2}.
Hence the substitution  $\Pi$  is an increasing substitution.
So for  $(l_{f-1}+2,t)$-fragment of the sequence $D_f$ we prove that
\begin{equation} 
\label{35} 
r_{f-1}=1, \, l_{f-1}=0, \, l_{f-2}>0. 
\end{equation}

Equality  (\ref{35}) with respect to the symmetry coinside with condition  (\ref{34}),
for a smalller value of $f$. 
By remark 3 we can apply the inductive assumption.
Lemma is proved.

Suppose that for a certain
 $i\geq0$ and for sequences  $\{b_j\}$ and
$\{c_j\}$ defined in  (4) one has
\begin{equation} 
\label{36} 
b_{i+1}=h_u\leq h_v=c_{i+1}
, 
\end{equation}
where  $\{h_i\}$ is the sequence from the definition of values  $W_f$. 
Put
$$H(i)=\{h_u,h_{u+1},\ldots,h_v\}
.
$$
Then $H(i)\neq\varnothing$.

\begin{Le}\label{l5}
Under the condition (\ref{36})  for the numbers 
$s_1,s_2,\ldots,s_k\in H(i)$,
$k\ge0$, one has
\begin{equation} 
\label{37-38} 
[0;n_i,n_{i-1},\ldots,n_0]-[0;s_1,s_2,\ldots,s_k,m_i,m_{i_1},\ldots,m_1]
\left\{
\begin{array}{ll}
\geq0 ,   \text{if } i\equiv0\pmod{2}, \\
\leq0 ,   \text{if } i\equiv1\pmod{2}, 
\end{array}
\right.
\end{equation}
\end{Le}

Proof.

By the definition for every $j=1,2,\ldots,k$
one has
$$
c_1\geq c_2\geq\ldots\geq c_{i+1}\geq s_j\geq b_{i+1}\geq
b_i\geq\ldots\geq b_1>b_0.
$$
That is why inequalities 
(\ref{37-38})
follows from the usisl rule of
comparing continued fractions' values step by step,
with respect to the parity of the step 
(of course here we should consider separate case of odd
 $k$ and even $k$). 
 Lemma is proved.

\subsection{Quadratic irrationalities' method.}
By induction define
$K_{l,n}$. 
Put
$$
K_{0,n}=1, \, K_{1,n}=n+2,\, K_{j+1,n}=(n+2)K_{j,n}-K_{j-1,n},\,\,
j=1,2,\ldots,n-1.
$$

\begin{Le}\label{l6}
For $l\geq1$ one has
\begin{equation} 
\label{41} 
\underbrace{\langle {1,n,1,n,\ldots,1}\rangle 
}\limits_{\mbox{$2l-1$ \text{elements}}}
=
K_{l-1,n}, 
\end{equation}
\begin{equation} 
\label{42} 
\underbrace{\langle {n,1,n,1,\ldots,n}\rangle 
}\limits_{\mbox{$2l-1$ \text{elements}}}
=
nK_{l-1,n}, 
\end{equation}
\begin{equation} 
\label{43} 
\underbrace{\langle 1,n,1,n,\ldots,1,n\rangle 
}\limits_{\mbox{$2l$ \text{elements}}}
=
K_{l,n}-K_{l-1,n}. 
\end{equation}
\end{Le}

This lemma obviosly follows by induction in  $l$.

Put
$$
k_{l,n}=\langle n^{\{ l\}}\rangle,\,\,\,\,\,\,
\lambda_n=\frac{n+\sqrt{n^2+4}}2
.
$$

\begin{Le}\label{l7}
The following equalities are valid:
\begin{equation} 
\label{44} 
k_{l,n}=
\left|\left|\frac{\lambda_n^{l+2}}{\lambda_n^2+1}\right|\right|, 
\,\,\,\,\,\,\,
K_{l,n}=\left[\frac{\mu_n^{l+2}}{\mu_n^2-1}
\right]
. 
\end{equation}
Here  $\mu_n$ is defined in  (\ref{22}).
\end{Le}

For  $K_{l,n}$
an equatity which was proved in \cite{7}
is very close to (\ref{44}).
We should note that our equality  (\ref{44}) for $k_{l,n}$ may be proved apsolutely analogouslely to those from \cite{7}.

\begin{Le}\label{l8}
For $n>8$ one has
\begin{equation} 
\label{45} 
\frac{1+(1/\lambda_n^2)}{1-(1/\mu_n^2)}\left(\frac{\mu_n}{\lambda_n}\right)^n<\frac{41}{40}e^2
.
\end{equation}
\end{Le}

To prove this lemma one should apply inequalities
$$
n\leq\lambda_n\leq{n+1}\leq\mu_n\leq{n+2}
$$ as well as the inequality $(1+1/n)^n<e$.

\begin{Le}\label{l9}
The following inequalities are valid:
\begin{equation} 
\label{46} 
K_{n,n}-K_{n-1,n}
<
k_{n+1,n}, \,\, K_{n-1,n}<k_{n,n}, 
\end{equation}
\begin{equation} 
\label{47} 
n>\frac{nK_{n-1,n}-k_{n,n}}{k_{n+1,n}-K_{n,n}}. \end{equation}
\end{Le}

Proof.
For  $n\leq8$
one can easily chech
(\ref{46}) and (\ref{47}).
For  $n>8$  
one should replace the values in the inequalities by the corresponding  values from Lemma \ref{l7}
and divide by
$\lambda_n^n$. Then one should apply (\ref{45}). Lemma is proved.

\subsection{
Majorizing inequalities' method.}

For $i,k\in\mathbb N$,
$j\in\mathbb N \cup \{0\}$  we consider
$$
L
=
\langle x_1,\ldots,x_i,z_1,\ldots,z_j\rangle 
,
$$
and
$$
M=
\langle y_1,y_2,\ldots,y_k,z_1,z_2,\ldots,z_j\rangle 
,$$
with the common tale 
 $(z_1,z_2,\ldots,z_j)$.

\begin{De}\label{d4}
We say that a continuant $M$ is majorized by a continuant $L$ (the notation $M\preceq L$)
if
\begin{equation} 
\label{48} 
\langle y_1,y_2,\ldots,y_k\rangle \leq \langle x_1,x_2,\ldots,x_i\rangle 
,
\end{equation}
\begin{equation} 
\label{49} 
\langle y_1,y_2,\ldots,y_{k-1}\rangle \leq
\langle x_1,x_2,\ldots,x_{i-1}\rangle  
.
\end{equation}
In the case when we have strict inequalities in  (\ref{48}), (\ref{49}) or the inequality in (\ref{48})
is a  strict one and 
$j\geq1$ we say that  $M$ is
strictly majorized by $L$ (notation
$M\prec L$).
\end{De}

\begin{Le}\label{l10}
For $i\geq2$ one has
$$
\langle x_1,x_2,\ldots,x_{i-1},1,z_1,z_2,\ldots,z_j\rangle \succ
\langle x_1,x_2,\ldots,x_{i-1}+1,z_1,z_2,\ldots,z_j\rangle 
.
$$
\end{Le}

A proof of lemma \ref{l10} easily follows from 
$$
\langle x_1,x_2,\ldots,x_{i-1},1
\rangle
=
\langle x_1,x_2,\ldots,x_{i-1}+1
\rangle 
,
$$
$$
\langle x_1,x_2,\ldots,x_{i-1}
\rangle 
>
\langle x_1,x_2,\ldots,x_{i-2}
\rangle 
.
$$

\begin{Le}\label{l11}
If $M\preceq L$ then $M\leq L$. 
Moreover if
$M\prec L$ then $M<L$.
\end{Le}

A proof of lemma \ref{l11} easily follows from 
the well-known formula
\begin{equation} 
\label{50} 
\langle x_1,\ldots,x_i,z_1,\ldots,z_j\rangle 
=
\langle x_1,\ldots,x_i\rangle 
\langle z_1,\ldots,z_j\rangle 
+
\langle x_1,\ldots,x_{i-1}\rangle 
\langle z_1,\ldots,z_{j-1}\rangle 
\end{equation}
(see \cite{11}) and from (\ref{48}), (\ref{49}).

\begin{Le}\label{l12}
Suppose that
 $jz_1>1$ and the inequality in  (\ref{48}) is a strict one. 
 In addition suppose that
\begin{equation} 
\label{51} 
\langle y_1,y_2,\ldots,y_{k-1},y_k+1\rangle 
\leq
\langle x_1,x_2,\ldots,x_{i-1},x_i+1\rangle 
. 
\end{equation}
Then $M\prec L$.
\end{Le}

Proof.
It is enough to show that
\begin{equation} 
\label{52} 
\langle y_1,y_2,\ldots,y_k,z_1\rangle \leq
\langle x_1,x_2,\ldots,x_i,z_1\rangle 
\end{equation}
holds and that this inequality is srtict for
$z_1>1$.

By (\ref{1}) we have
$$
\langle y_1,y_2,\ldots,y_k,z_1\rangle 
=
\langle y_1,y_2,\ldots,y_k,1\rangle +(z_1-1)
\langle y_1,y_2,\ldots,y_k\rangle 
,
$$
$$
\langle x_1,x_2,\ldots,x_i,z_1\rangle 
=
\langle x_1,x_2,\ldots,x_i,1\rangle +
(z_1-1)\langle x_1,x_2,\ldots,x_i\rangle 
.
$$
By  (\ref{32}) and (\ref{51})
we compare the corresponding therm in the two last inequalities and obtain
(\ref{52}). 
Lemma is proved.

\begin{Le}\label{l13}
Suppose that we have a strict inequality in 
(\ref{48}). 
Moreover supppose that
$$
[z_j;z_{j-1},z_{j-2},\ldots,z_2,z_1]
>
\frac{\langle y_2,\ldots,y_k\rangle 
-
\langle x_2,\ldots,x_i\rangle }{\langle x_1,\ldots,x_i\rangle -\langle y_1,\ldots,y_k\rangle }
.
$$
Then
\begin{equation} 
\label{53} 
\langle z_1,z_2,\ldots,z_j,x_1,x_2,\ldots
x_i\rangle 
>
\langle z_1,z_2,\ldots,z_j,y_1,y_2,\ldots,y_k
\rangle 
. 
\end{equation}
\end{Le}

The proof easily follows by applying equalities (\ref{50})
to continuants from (\ref{53}) and by taking into account (\ref{2}).

Put
$$
J(l)=\{\underbrace{1,n,1,n,\dots,1}_{(2l-1)\, elements}\},\,\,\,
G=\left\{z,n,J(m_-)\right\}
.
$$

\begin{Le}\label{l14}
For every $n,m\in\mathbb N$ one has
\begin{equation} 
\label{54} 
[n;J(m_+)]
>
\frac{\langle
n,J(n-1),n\rangle-\langle
n^{
\{n\}
}\rangle}{\langle n^{
\{n+1\}
}\rangle-\langle J(n),n\rangle}
,  
\end{equation}
\begin{equation} 
\label{55} 
\langle
n^{
\{n\}
}\rangle
>\langle J(n)\rangle,
\end{equation}
\begin{equation} 
\label{56} 
\langle  n^{
\{n+1\}}\rangle
>
\langle J(n),n\rangle
\end{equation}
\begin{equation} 
\label{57} 
\langle  J(m_+),
n^{
\{n+1\}}\rangle
>
\langle  J(m_+),n,J(n)\rangle
\end{equation}
\begin{equation} 
\label{58} 
\langle  J(m_+),  n^{
\{n+2\}}\rangle
>
\langle  J(m_+),n,J(n),n\rangle
.
\end{equation}
Moreover for
$x\geq n+1$ one has
\begin{equation} 
\label{59} 
\langle  J(m_+),  n^{
\{x+1\}},G\rangle
>
\langle  J(m_+),n,J(n),  n^{
\{x-n\}},G\rangle
. 
\end{equation}
\end{Le}

Proof.

Inequality  (\ref{54}) follows immediately from (\ref{47}) and Lemma  \ref{l6}.
Inequalities (\ref{55}) and (\ref{56}) are easy corollaries of
(\ref{41}), (\ref{43}) and (\ref{46}).
This proves the majorizind inequality
\begin{equation} 
\label{60} 
\langle  n^{\{ n+1\}},J(m_+)
\rangle
\succ
\langle  J(n),n,J(m_+)
\rangle
. 
\end{equation}
By the symmetry and Lemma \ref{l11} we see that  (\ref{60}) 
leads to (\ref{57}).
To prove  (\ref{58})
one should apply Lemma \ref{l13}. 
Conditions od Lemma \ref{l13} are satisfied as we have (\ref{54}) and (\ref{56}).
Then inequalities  (\ref{57}) and  (\ref{58})
lead to strict majorizing property in the inequality (\ref{59}) which now follows from Lemma \ref{l11}.
Lemma is proved.

\subsection{ Unit variation method.}
The idea of this method is taken from \cite{7}.

\begin{Le}\label{l15}
Consider the sequence
$Q$ which gives the minimal value of the continuant on a set of the form $U(S,n,t)$
or $U_n(S)$.
Then in this sequence there exist not more than one element different from 
$1$ and 
$n$.
\end{Le}

Proof.
Assume that it is not so.
Then we can take two elements $a$, $b$ from  $Q$
which are not equal to $1$ or $n$.
 Now we suppose all other elements to be fixed.
We  suppose the sum 
$a+b=\tau$ to be fixed also. 
Consider the continuant
$$
F_\tau(a)
=
\underbrace{\langle  \ldots,a,\ldots,\tau-a,\ldots\rangle}_{t\,\, elements}  
.
$$
Then from  (\ref{50}) we see that the function $F_\tau(a)$
is a quadratic polynomial in $a$ with a negative leading coefficient.
Hence by convexity argument
 $$
F_\tau(a)>min\left\{F_\tau(a-1),\,\,\,F_\tau(a+1)\right\}
$$
(here one can see an explanation for the notion ``unit variation``).
This inequality contradicts the minimality of $Q$. Lemma is proved.

\begin{Le}\label{l16}
Suppose that $a,b,c\in\mathbb N$, $p\in\mathbb N\cup\{0\}$. 
If
$a\geq c+1$ then
\begin{equation} 
\label{61} 
\langle a,b^{\{ p\}},c\rangle
\leq 
\langle
a-1,b^{\{ p\}},c+1\rangle
. 
\end{equation}
If in addition $a>c+1$ then there is a strrict inequality in (\ref{61}).
\end{Le}

Proof.

When
 $a=c+1$ the equality  (\ref{61})
is obvious as we have symmetry property (\ref{31}).
In the case  $a>c+1$ by  (\ref{32}) and Lemma \ref{1A} we have
$$
\langle a,b^{\{ p\}},c\rangle
=
\langle 1,a-1,b^{\{ p\}},c\rangle
<
\langle 1,c,b^{\{ p\}},a-1\rangle
=
\langle
a-1,b^{\{ p\}},c+1\rangle
.
$$

\begin{Le}\label{l17}
If  $a\geq c+2$ and $p\geq0$, $jz_1>1$ then
$$
\langle a,b^{\{ p\}},c,z_1,z_2,\ldots,z_j\rangle
\prec
\langle a-1,b^{\{ p\}},c+1,z_1,z_2,\ldots,z_j\rangle
.
$$
\end{Le}

To proof this lemma one should apply Lemma \ref{l12}.
To do this we must verify the assumotion hypothesis of this lemma.
This cam be done by Lemma \ref{l16}.

\subsection{ Lemmas concerning systems of (in)equalities.}

\begin{Le}\label{l18}
Consider $S,t\in\mathbb N$. 
The system  
\begin{equation} 
\label{vi9o13} 
\left\{
\begin{array}{ll}
    ch_1+dh_2=S,\\
    c+d=t, \\
    h_2=h_1+1,\\
     d<t
\end{array}
\right.
\end{equation} 
has the unique solution in integers 
$c,d,h_1,h_2$.
If in addition $S\not \equiv0  \pmod t $ then
this solution is defined by equalities  (\ref{15});
in the opposite case
\begin{equation} 
\label{56c4615} 
h_1=\frac{S}{t}, \, \,\,
h_2=h_1+1, \, c=t , \, d=0
. 
\end{equation}
\end{Le}

Proof. 
According to the equations of the system (\ref{vi9o13}), 
$$
S= 
ch_1+dh_2
=
ch_1+d(h_1+1)
=
(c+d)h_1+d
=
th_1+d
.
$$
Therefore, $d$ is the remainder of dividing $S$ by
$t$, while $h_1$ is the quotient of this division.
Hence the formulas (\ref{15}) and (\ref{56c4615}). 
The lemma is proved.

\begin{Le}\label{l19}
Suppose that $2\leq t\leq S<nt$. 
Then the only one from the following two systems 
\begin{equation} 
\label{g6c4615} 
\begin{array} {lr}
\left\{
\begin{array} {l}
m(n+1)+nx+z=S,\\
2m+x+1=t,\\
m\geq0,\\
x\geq0,\\
1\leq z\leq n-1,
\end{array}
\right. &
\left\{
\begin{array}{l}
m(n+1)-x+z=S,\\
2m-x+1=t,\\
m\geq0,\\
x<0,\\
1\leq z\leq n-1.
\end{array}
\right.
\end{array}
\end{equation}
has a solution in integers $m,x,z$. 
This solution is unique. 
It is determined be formulas (\ref{17}), (\ref{18}).
\end{Le}

Proof.
The second equality in (\ref{18}) is obtained from the second equalities of each of the two systems 
in  (\ref{g6c4615}).
It also shows that the parity of the numbers $t$ and $x$ are opposite.
When substituting the second equality in (\ref{18}), taking into account the sign of the variable $x$, in the first equations of each of the two systems in  (\ref{g6c4615}). 
Then the equality is obtained
\begin{equation}
\label{19c46189}
z=
S -t+1
-\frac{t-1+x}{2}(n-1)
.
\end{equation}
According to the above, the numerator of the fraction in
(\ref{19c46189}) is even, so that the coefficient at $(n-1)$  is whole.
From here we get a comparison in
(\ref{17}).

We show that the number $m$ satisfying the equalities in any of the systems in  (\ref{g6c4615}) is non-negative.
To beging with, let us  consider the first of these two  systems.
Expressing $x$ from the second equality and substituting it into the first one, we get
$$
m(n+1)+n(t-1-2m)+z
=S
<nt
.
$$
Simplifying the resulting inequality, we derive:
$$
m>\frac{z-n}{n-1}\geq -1
.
$$

Similarly, expressing $x$ from the second equality of the second system  in  (\ref{g6c4615}) and substituting  it into the first of them, we get
$$
m(n+1)+(t-1-2m)+z
=S
\geq t
.
$$
Simplifying the resulting inequality, we derive:
$$
m\geq \frac{1-z}{n-1}> -1
.
$$
Hence, since $m$ is an integer, we get the inequality $m\geq 0$.

Finally, expressing $x$ of (\ref{19c46189}) obtained first equality in (\ref{18}).
Depending on the sign of $x$
triple $(m,x,z)$, found by far,
be  solution of the first or the second systems   in  (\ref{g6c4615}).
The lemma is proved.

For a natural $T$, we put 
\begin{equation} 
\label{19c4674} 
T_0^{(n)}
=
(n+1)\left\{
\frac{T-1}{n+1}
\right\}
,
\end{equation}
\begin{equation} 
\label{19c46744} 
T_1^{(n)}
=
n\left\{
\frac{T}{n}
\right\}
.
\end{equation}

\begin{Le}
\label{l29o51} 
For
$T<n^2-n-1$, the inequality
$
T_1^{(n)}
\leq 
T_0^{(n)}
$
is satisfied if and only if
the relation holds
\begin{equation}
\label{1046744}
\left[
\frac{T}{n}
\right]
-
\left[
\frac{T-1}{n+1}
\right]
=1  
\end{equation}
(otherwise the left side of the equation 
(\ref{1046744}) is zero).
\end{Le}

Proof. 
Let $T_1^{(n)}\le T_0^{(n)}$. 
We prove the equality (\ref{1046744}). 
On the one hand,
\begin{equation} 
\label{11vg4} 
\begin{split} 
n
\left[
\frac{T}{n}
\right]
=
T-
n\left\{
\frac{T}{n}
\right\}
=
T-
T_1^{(n)}
\geq
T-
T_0^{(n)}
=
T-
(n+1)\left\{\frac{T-1}{n+1}\right\}
=
\\=
(n+1)\left[\frac{T-1}{n+1}\right]
+1
,
\end{split} 
\end{equation}
from where, dividing the beginning and end of the chain of inequalities (\ref{11vg4}) by $n$, we get
\begin{equation} 
\label{1146g4} 
\left[
\frac{T}{n}
\right]
\geq
\frac{n+1}{n}
\left[\frac{T-1}{n+1}\right]
+\frac{1}{n}
>
\left[\frac{T-1}{n+1}\right]
.
\end{equation}
On the other hand,  
\begin{equation} 
\label{104604} 
\left[
\frac{T}{n}
\right]
-
\left[
\frac{T-1}{n+1}
\right]
=
\frac{T}{n}
-
\frac{T-1}{n+1}
-
\left\{
\frac{T}{n}
\right\}
+
\left\{
\frac{T-1}{n+1}
\right\}
=
\frac{T+n}{n(n+1)}
-
\left\{
\frac{T}{n}
\right\}
+
\left\{
\frac{T-1}{n+1}
\right\}
.
\end{equation}
And since $T<n^2$, then (\ref{104604}) leads to the inequality$$
\left[
\frac{T}{n}
\right]
-
\left[
\frac{T-1}{n+1}
\right]
\leq
\frac{T+n}{n(n+1)}
+
1
<
2
.
$$
From here and from (\ref{1146g4}) it follows (\ref{1046744}).

Now let the inequality $T_1^{(n) }> T_0^{(n)}$ be satisfied, or
$$
n\left\{\frac{T}{n}\right\}
> 
(n+1)\left\{
\frac{T-1}{n+1}
\right\}
.
$$
Dividing both parts of the last inequality by $n$, we get:
$$
\left\{\frac{T}{n}\right\}
> 
\frac{n+1}{n}
\left\{\frac{T-1}{n+1}\right\}
\geq
\left\{\frac{T-1}{n+1}\right\}
.
$$
From here and from (\ref{104604}) for $0< T<n^2$ we get
$$
-1
<
\left[
\frac{T}{n}
\right]
-
\left[
\frac{T-1}{n+1}
\right]
<
1
,
$$
from where it follows
$$
\left[
\frac{T}{n}
\right]
-
\left[
\frac{T-1}{n+1}
\right]
=0
.
$$
The lemma is proved.

\begin{Le}
\label{l29o11} 
First,
every integer 
$T$ greater than $n^2-n-1$ is represented as the sum of a certain number of terms that can take the values $n$ or $n+1$; the number $T=n^2-n-1$ does not have this property. 

Second, the solution of the equation $T=m(n+1)+nx$ in the non-negative integer variables $m$ and $x$ exists if and only if its solution exists with an additional constraint 
$x\leq n$. 

Third,
for $T< n^2-n-1$, the number $T$ is represented as the sum of a certain number of terms
of the form $n$ or $n+1$ if and only if 
\begin{equation} 
\label{1146744} 
T=
(n-T_0^{(n)})n
+
T_1^{(n)}(n+1)
.
\end{equation}

Fourth, the equality (\ref{1146744}) holds
if and only if
$\,\,\,T_1^{(n)}\leq T_0^{(n)}$. 
\end{Le}

Proof. 
The first part of the Lemma is a particular case of approval of Sylvester \cite{Syl}: 
every integer greater than $ab-a-b$, can be represented as the sum of a certain number of terms, which can take the values $a$ or $b$ (it is sufficient to substitute $a=n$, $b=n+1$), while the number of $ab-a-b$ is not representable in this form.

We prove the second part of the lemma. 
To do this, divide $x$ with the remainder by $n+1$:
let
$$
x=q(n+1) +r,\,\,\,\,r\leq n. 
$$
When substituting the last equality in the Diophantine equation, we get:
$$
T=
m(n+1)+nx
=
m(n+1)+n(q(n+1) +r)
=
(m+nq)(n+1)+nr
,
$$
and we got a solution with the desired property.

We prove the third part of the lemma. 
We show that if
the number $T<n^2$ is represented by the sum of terms of the form $n$ or $n+1$, then
the equality (\ref{1146744}) is satisfied. 
Indeed: let 
\begin{equation} 
\label{vi9os1} 
T=m(n+1)+nx
\end{equation}
and $m\geq 0$, $x\geq 0$.
Then, taking the remainder of the division of the left and right parts of the equality (\ref{vi9os1}) by $n$, we get: $m\equiv T \pmod {n}$, which implies that $m=T_1^{(n)}$. 
Similarly, taking the remainder of the division of the left and right parts of the equality (\ref{vi9os1}) by $n+1$, we get: $x\equiv -T \pmod {n+1}$, which implies that $x=n-T_0^{(n)}$. 
Substituting the found values of $m$ and $x$
in (\ref{vi9os1}) results in equality (\ref{1146744}).

We prove the fourth part of the lemma. 
We now show that for $T<n^2-n-1$
, the number $T$ is represented by the sum of terms of the form $n$ or $n+1$
if and only if $T_1^{(n)}\leq T_0^{(n)}$.
To do this, transform the right side of the equality (\ref{1146744}):
$$
(n-T^{(n)}_0)n+T^{(n)}_1(n+1)
=
\left(n-
(n+1)\left\{\frac{T-1}{n+1}\right\}
\right)n
+
n\left\{\frac{T}{n}\right\}(n+1)
=
$$
$$
=
n^2
-(n^2+n)\left(
\left\{\frac{T-1}{n+1}\right\}
-\left\{\frac{T}{n}\right\} 
\right)
=
$$
$$
=
n^2
-(n^2+n)
\left(
\frac{T-1}{n+1}
-
\left[\frac{T-1}{n+1}\right]
-\frac{T}{n}
+\left[\frac{T}{n}\right] 
\right)
=
$$
$$
=
n^2-n(T-1)+(n+1)T
-(n^2+n)
\left(
\left[\frac{T}{n}\right] 
-
\left[\frac{T-1}{n+1}\right]
\right)
=
$$
$$
=
T
+
(n^2+n)
\left(
\left[\frac{T}{n}\right] 
-
\left[\frac{T-1}{n+1}\right]
-1
\right)
.
$$
Thus, the right side of an equality (\ref{1146744}) is equal to the left one if and only if the equality (\ref{1046744}) is satisfied, which, according to the previous lemma, is satisfied if and only if $T_1^{(n)}\leq T_0^{(n)}$. 
The lemma is proved.

\begin{Le}\label{l20}
Suppose that  $2\leq t\leq S<nt$.
Then  for any 
 $z\in P(S)$
there exists only one solution of the system
\begin{equation} 
\label{62} 
\left\{
\begin{array}{l}
m(n+1)+nx+z=S,\\
m\geq0\\
 0\leq x\leq n
\end{array}
\right. 
\end{equation}
in integers $m,x$. 
This solution is defined by formulas (\ref{20}).

Moreove an opposite statement is true:
if $S\leq n^2-1$, $z\in\{1,2,\ldots,n-1\}$ and the system
(\ref{62}) has a unique solution then $z\in P(S)$.
\end{Le}

Proof. 
The number $T=S-z$, according to the first equation of the system (\ref{62}),
is represented as the sum of a certain number of terms of the form $n$ or $n+1$. 
Therefore, according to the preceding lemma, the inequalities $S-z>n^2-n-1$ or
$$
n
\left\{\frac{S-z}{n}\right\}
\leq
(n+1)
\left\{\frac{S-z-1}{n+1}\right\}
$$
are satisfied.
The first of these corresponds to the case when $z\in P(S)$ as first line 
in (\ref{8}); 
the second and third lines in (\ref{8}) correspond to the situations when $S<n^2-1$, the number $S$ is not representable or representable 
as the sum of several terms, equal to $n$ or $n+1$. 
The uniqueness of the solution is provided by the last inequality in (\ref{62}). 

The expressions for $x$ and $m$ in (\ref{20})
follow from the first line of the system (\ref{62}). 

The lemma is proved.

\section{ Proof of main theorems.}

\subsection{Proof of Theorem 1. }
Induction in  $t$.
Suppose
 $A=\{a_1,a_2,\ldots,t\}$ to be a sequence from $V_f$.
Suppoce than any basic substitution for it is a non-increasing substitution.
Then $a_1=b_0$  by the assumption. If $a_2=b_1$ then $(2,t)$-fragment of the sequence 
$A$ sshould belong to a certain set  $V_g$ for a certain value of $g\leq f$.
Moreover it containes of  $t-1$ elements.
We apply the inductive assumotion to this fragment.
We  should take into account remark 3 also.
We see that  the equality  (\ref{9})
is proved, in this case.
If $a_2\neq b_1$ (and hence 
$a_2>b_1$) we consider the basic  $(2,\gamma)$-substitution
$\Pi$, $\gamma=\max\{j:a_j=b_1\}$.
Then $\Pi$ is an increasing substitution  by Lemma \ref{1A}. 
This contradicts to the assumption behind.
So $a_2=b_1$. 
Theorem is proved.

\subsection{Other proof  of Theorem \ref{t1}. }
Consider three finite sequences consisting of natural numbers: 
$U,V$, and $W$, the middle of which (i.e., $ V$) is --- nonempty, as is at least one of the two extremes, --- 
$U$ or $W$. 
Let us give the following definition.

\begin{De}
\label{3e2sf}
Let's call the procedure for replacing continuants 
$$
\langle { U,V,W} \rangle 
 \,\,\,\, \mapsto\,\,\,   
\langle { U},\overleftarrow{V} , {  W}  \rangle
$$
{\it by reflection} \cite{Gay} of the sequence $   ( U,V,W)$, where $\overleftarrow{V}$ --- the sequence $V$ written in the opposite order.
\end{De}

Consider any finite sequence
$(a_1,a_2\dots,a_t)$ of natural numbers, such that 
$ a_1=\min(a_1,a_2,$ $\dots,a_t)$ and 
$ a_i >a_j$ for some $i<j$ (any such pair $(a_i,a_j)$ let's call {\it correct}). 
Let $ a_{i^*}, a_{j^*}$ --- two elements of this sequence, such that 
$ a_{i^*} >a_{j^*}$, ${i^*}<{j^*}$, with the maximum difference value
$ j^* - i^*$ (and such a correct pair $(a_{i^*}, a_{j^*})$ let's call {\it the most remote}). 
Then $i^*>1$ (since $ a_1$ is minimal) and
\begin{equation}
\label{tuda}
a_{i^*-1} \le a_{j^*} <a_{i^*} \le a_{j^*+1} 
\end{equation}
(of course, it may happen that $j^*=t$;
then
the sequence $W$ defined below is empty).
  
Let
\begin{equation}
\label{tu4da}
\langle  a_1, a_2,\dots,a_t\rangle 
= 
\langle U,V,W\rangle
,
\end{equation}
where 
\begin{equation}
\label{tu3da}
U = a_1, a_2,\dots,a_{i^*-1}
;
\,\,\, 
V= a_{i^*},\dots,a_{j^*} 
;
\,\,\,\, 
W= a_{j^*+1},\dots,a_t
.
\end{equation}
It follows from (\ref{tuda}) that the inequality 
\begin{equation}
\label{tu1da}
{\bf a}(U,V,W)
\ge
0
\end{equation}
is satisfied.
Therefore, it follows from the Lemma \ref{1B} that 
\begin{equation}
\label{tu2da}
\langle { U},\overleftarrow{V} , {  W}  \rangle
\geq 
\langle
U,V,W\rangle
.
\end{equation}
However, the sequence 
$ ({ U},\overleftarrow{V} , {  W})$ can be obtained from the sequence $(a_1, a_2\dots, a_t)$
by some permutation (reflection).
Thus, for each sequence
$(a_1,a_2\dots,a_t)$ with inversion $ a_i >a_j, i<j$ 
there is a permutation $\pi$ that reduces the number of inversions, but
does not reduce the continuant $\langle a_1,a_2,\dots,a_t \rangle$.
The theorem \ref{t1} is proved.

In general terms, this proof consisted in the fact that from an arbitrary continuant, using an algorithm consisting of non-decreasing reflections, a continuant equal to the maximum over all permutations of its elements was obtained. 
For the rest, let's imagine the following situation: let the continuant (\ref{tu4da}), which we will call the {\it internal} continuant, be part of another, --- {\it external}, --- continuant (say, $\langle
F, a_1, a_2,\dots,a_t, G \rangle
$), where $F$ and $G$ are finite sequences. 
Suppose that for some purposes it is required
to obtain a similar algorithm, which from an arbitrary continuant leads to a maximum, but consists
only of such reflections that
do not decrease for both the internal continuant and the external one. 
This algorithm is called {\it transitive} (step-up). 
Similarly, a transitive step-down algorithm is defined: it consists of reflections that are not external  continuant, as for the internal one.

\begin{Le}\label{ltud7a}
For any continuant
that is internal
with respect to some
external one, a transitive algorithm exists (both increasing and decreasing).
\end{Le}

Proof. 
We prove, for example, the existence of a transitive step-up algorithm. 
Namely, if for the next step of the proof the constructed reflection (\ref{tu2da}) strictly increases the internal continuant, then it also strictly increases the external one. 
Indeed: from the rule of comparing continued fractions and from the positivity of a number 
$ 
{\bf a}(U,V,W) 
$ from 
(\ref{26b})
follows also the positivity of the number 
$ 
{\bf a} ((F,U),V,(W,G)) 
$.
Therefore, it is sufficient to deal only with those reflections that do not change the value of the internal continuant: in this case, the inequalities in (\ref{tu1da})
and (\ref{tu2da}) are performed as equalities. 
Thus, 
${\bf a} (U,V,W) $ from 
(\ref{tu1da}) is zero. 
But $ {\bf a} (U,V,W) $ is equal to the product of two multipliers. 
Hence, the general situation falls into two cases:
$
[\overleftarrow{ V}]=[{V}]
$
or
$[\overleftarrow{U}]=[{W}]
$.
However, according to (\ref{tuda}), the first of these cases is not possible.

Let it be now $
[\overleftarrow{ V}]\not=[{V}]
$, but
the equality is met
$[\overleftarrow{ U}]=[{{ W}}]
$.
According to (\ref{tuda}), the final sequences $\overleftarrow{ U}$ and $W$
differ already by their first elements, so
that ``unit extraction"  (\ref{32}) is involved here.
More precisely, $\overleftarrow{ U}=\{n, 1\}$, $W=\{n+1\}$ (a singleton sequence) for some natural $n$.
We show that in this case, the following is true:
\begin{equation}
\label{tud7a}
a_{i^*-1} = a_{j^*}
<
a_{i^*} = a_{j^*+1}
.
\end{equation}
Indeed: if $a_{i^*}>a_{j^*+1}$, then
the correct maximally distant pair --- is
$(a_{i^*}, a_{j^*+1})$ , not
$(a_{i^*} , a_{j^*})$.
If $a_{i^*}<a_{j^*+1}$, then
the correct maximally distant pair is
$(a_{i^*-1} ,a_{j^*})$ , not
$(a_{i^*}, a_{j^*})$.
So $a_{i^*}=a_{j^*+1}=n+1$, while $a_{i^*-1}=n$.
So $a_{j^*}<n+1$.
But, if $a_{j^*}<n=a_{i^*-1}$, then, again,
the correct maximally distant pair --- is
$(a_{i^*-1} ,a_{j^*})$ , not
$(a_{i^*}, a_{j^*})$.
This contradiction finally proves the equalities in (\ref{tud7a}).

Now $[V]<
[\overleftarrow{ V}]$ (due to the inequality between the first and last elements of the final sequence $V$), and if at least one of the two sequences $F$ or $G$ is non-empty, then
$[\overleftarrow{U},\overleftarrow{F}]>
[W,G]$.
Therefore, for the outer continuant, the reflection corresponding to the reflection of the inner one (\ref{tu2da}) is strictly magnifying.
Thus, the existence of a transitive step-up algorithm is proved.

The existence of a transitive step-down algorithm is proved in a completely analogous way.

The lemma is proved.

\subsection{Proof  of Theorem \ref{t2}. }

By Lemma \ref{l3} 
the maximum under consideration attains on a certain sequence $D_f$.
By means of the symetry one may assume that at least one of two numbers $l_f$, $r_f$,
(say, the second one) is greater than zero.
In order to obtain a proof of Theorem \ref{t2} from  Lemma \ref{l4}
we add the continuant by a virtual element
$h_{f+1}=+\infty$, formally:
$$
\langle+\infty,D_f\rangle
=
\langle D_f\rangle
.
$$
Then we have
$$
l_{f+1}=1,\,\,\, r_{f+1}=0,\,\,\, r_f>0.
$$
This is just  (\ref{34}) for large  $f$. 
Now Theorem \ref{t2} folows from Lemma \ref{l4}.
Theorem is proved.

\subsection{Other proof  of Theorem \ref{t2}. }

\begin{Le}
\label{l3} 
Maximal continuant on the set 
$W_f(\bar{h},\bar{p})$
is reached on some finite sequence 
$D_f(\bar{l},\bar{r})$, whose elements are subject to the condition (\ref{3}).
\end{Le}

Proof. 
Let be the maximal continuant on the set 
$W_f(\bar{h},\bar{p})$
is reached on a finite sequence
$D=(a_1,a_2,\ldots,a_t)$ and let the equality $a_\nu=h_1$ be satisfied for some $\nu\in\{1,2,
\ldots,t\}$. 
Let's imagine the sequence $D$ as a
union of finite subsequences 
$(a_1,a_2,$ $\ldots,a_\nu)$ and $(a_\nu,a_{\nu+1},\ldots, a_\nu)$, in the first of which the minimum element is in the last place, and in the second --- in the
first place. 
We fix the location of the minimum element in each of these sequences and apply to them separately the transitive algorithm for proving the \ref{t1} theorem (in the first of them, changing the order of the elements to the exact opposite, and at the end of the algorithm --- reversing it again). 
Then the first sequence turns into a non-increasing one, the second --- into a non-decreasing one, and together they form a finite sequence 
$D_f(\bar{l},\bar{r})$.
The lemma is proved.

Consider an arbitrary continuant of the form
\begin{equation}
\label{19c46350}
\langle D_f\rangle
=
\langle
h_f^{\{l_f\}},h_{f-1}^{\{l_{f-1}\}},\ldots,h_2^{\{l_2\}},
h_1^{\{l_1\}},h_1^{\{r_1\}},h_2^{\{r_2\}},\ldots,h_f^{\{r_f\}}
\rangle
\end{equation}
together with a symmetric one to it, which is obtained by writing all its elements in exactly the opposite order.
From these two continuants, we choose the one for which $l_{f-1}\geq r_{f-1}$.
In the continuant (\ref{19c46350}), we select the inner part:
\begin{equation}
\label{19c46351}
\langle D_{f-1}\rangle
=
\langle
h_{f-1}^{\{l_{f-1}\}},h_{f-2}^{\{l_{f-2}\}},\ldots,h_2^{\{l_2\}},
h_1^{\{l_1\}},h_1^{\{r_1\}},h_2^{\{r_2\}},\ldots,h_{f-1}^{\{r_{f-1}\}}
\rangle
\end{equation}
and apply a transitive step-up algorithm to it, using the lemma \ref{ltud7a}.
This means that by applying a certain number of reflections
, the continuant (\ref{19c46351}) will be brought to the maximum over all permutations of its elements, while the continuant (\ref{19c46350}) will not decrease.
Without limiting generality, we can assume that the inequality $l_{f-1}\geq r_{f-1}$ still holds.

Using the inductive (on the parameter $f$) assumption, we obtain that for the continuant (\ref{19c46351}) at least one of the sets of equalities is satisfied, 
(\ref{11}) or the opposite to it --- the one for which the inequality $l_{f-1}\geq r_{f-1}$is satisfied. 

If now $l_f=1$, then the inductive transition is complete. 
It remains to consider the cases when $l_f=0$ or $l_f>1$
Consider the case where $l_f=0$. 
Let's put it:
$
U$
--- 
empty,
$$
V
=
\left\{  
h_{f-1}^{\{l_{f-1}\}},h_{f-2}^{\{l_{f-2}\}},\ldots,h_2^{\{l_2\}},
h_1^{\{l_1\}},h_1^{\{r_1\}},h_2^{\{r_2\}},\ldots,
h_{f-1}^{\{r_{f-1}\}}
,
h_{f}^{\{r_{f}-1
\}}\right\},\,\,\,
W=\left\{ h_{f}\right\}.$$
Applying a single reflection with such a partitioning of the sequence into $U,V, W$ increases the continuant, so that now $l_f>0$. 

Now let $l_f>1$. 
Let's put it:
$$
U=\left\{  
h_{f}^{\{l_{f}-1\}}
\right\},
$$
$$
V=\left\{h_{f},  
h_{f-1}^{\{l_{f-1}\}},h_{f-2}^{\{l_{f-2}\}},\ldots,h_2^{\{l_2\}},
h_1^{\{l_1\}},h_1^{\{r_1\}},h_2^{\{r_2\}},\ldots,
h_{f-1}^{\{r_{f-1}\}}
\right\},\,\,\,$$ $$
W=\left\{ h_{f}^{\{r_{f}
\}}\right\}
.
$$
Applying a single reflection with such a partition of the sequence into $U, V, W$ increases the continuant, so that now $l_f=1$. 
The inductive step is complete. 
The \ref{t2} theorem is proved.

\subsection{Proof  of Theorem \ref{t3}. }
Assume that theorem \ref{t3} is not true.
So the minimal value of the continuant on the set $W_f$
is attained on a sequence
$B=\{a_1,a_2,\ldots,a_t\}$,
and no one of the trivial substitutions od this sequence does not give  (\ref{12}).
Then there exists
$$
j=j(B)=\min\{k\geq0:(a_{k_-+1}\neq
n_{k_-},k_-\leq\nu)\bigvee(a_{k_-+1}=n_{k_-},a_{t+1-k_t}\neq
m_{k_+},1\leq k_+\leq\mu)\}.
$$
Here as usually  $k_-=\left[  k/2\right]$, $k_+=k-k_-$.
We consider all the sequences 
trivially obtained from $B$.
We may choose $B$ to be the sequence with the maximal value of  $j(B)$.
Two cases are possible:
  $a_{j_-+1}\neq
n_{j_-}$ or $(a_{j_-+1}=n_{j_-},a_{t+1-j_+}\neq m_{j_+})$. 
In the first case  by $\nu$
we denote the index for which $a_\nu=n_{j_-}$.
Then we consider
the basic $(j_-+1,\nu)$-substitution  $\Pi$. 
In the second case  we put $a_\nu=m_{j_+}$. In this case we take   
 $\Pi$ to be the  $(\nu,t+1-j_+)$-basic substitution.
By Lemmas  \ref{1B},   \ref{l5} we see that the substitution  $\Pi$ is a non-increasing substitution.
From the minimality property we see that it is a non-decreasing
substitution also.
Hence by  Remark  \ref{z2} 
$\Pi$ is a trivial substitution.
But by the construction we have
$$
j(\{a_{\Pi(1)},a_{\Pi(2)},\ldots,a_{\Pi(t)}\})>j(B).
$$
This contradicts to the choise of 
$j$. 
Theorem is proved.

\subsection{Other proof  of Theorem \ref{t3}. }

We will conduct induction along the length of the continuant.
To do this, first consider the element $a_1$, selecting the desired numbering of the elements --- from the first to the last or from the last to the first, so that the inequality is satisfied
$a_1\leq a_t$.
If $a_1$ --- is not the minimum element of the continuant, then we apply a downward reflection that moves this minimum element to the first place.
Let this reflection result in the continuant $\langle a '_1,a' _2,\ldots,a '_t\rangle $, where $a'_1=h_1$.
If now $a'_2\not=h_f$, that is, to the maximum element of the continuant, then we apply a lowering reflection that moves this maximum element to the second place.
Let this action result in the continuant $\langle 
a'' _1,a''_2,\ldots,a''_t\rangle $.
Note that the elements $a''_1$ and $a'' _2$ coincide, respectively, with $n_0$ and $n_1$ --- 
the first two elements of the continuant from the right side in (\ref{12}).

Consider a finite sequence $ a'' _1,a''_2,\ldots,a''_t$ and change it numeration reverse (from end to beginning), applying the inductive hypothesis in the form of a transitive reduction algorithm, then you come to formula (\ref{12}).
Theorem \ref{t3} is proved.

\subsection{Proof  of Theorem \ref{t4}. }
By Lemmas \ref{l10} and \ref{l11}
applied to continuants
corresponding to  sequences containing  different from $ 1$  elements we see that
the equality (\ref{9}) should be valid. 
Theorem is proved.

\subsection{Proof  of Theorem \ref{t5}. }

Let the finite sequence $C$ for which
the maximum of continuants is reached on the set $U(S,t)$
belong to the set
$W_f$ with $f\geq3$.
Then the structure of this sequence is given by the formulas (\ref{10}) --- (\ref{11}).
But, according to the Lemma \ref{l17},
$$
\langle h_f,h_{f-1}^{\{ l_{f-1}\}},h_{f-2},\ldots\rangle
\prec
\langle h_f-1,h_{f-1}^{\{ l_{f-1}\}},h_{f-2}+1,\ldots\rangle.
$$
The last inequality, according
to the Lemma \ref{l11}, comes into conflict with the maximality of the continuant on the sequence $C=D_f$.
This contradiction shows that the inequality must hold
$f\leq 2$.

We show that in the case of $f=2$, the equality $h_2=h_1+1$ must hold.
Indeed, if $h_2>h_1+1$, then
\begin{equation}
\label{vi9oh1}
\langle h_2,h_{1}^{\{ l_{ 1}\}}
,h_2,\ldots\rangle
\prec
\langle h_2,h_{1}^{\{ l_{ 1}-1\}},
h_{1}+1,h_2-1,
\ldots\rangle
,
\end{equation}
since
$$
\langle h_2,h_{1}^{\{ l_{ 1}\}}
\rangle
<
\langle h_2, h_{1}^{\{ l_{ 1}-1\}},
h_{1}+1\rangle
$$
and (from (\ref{50}))
$$
\langle h_2,h_{1}^{\{ l_{ 1}\}}
,h_2 \rangle
=
\langle h_2,h_{1}^{\{ l_{ 1}-1\}}\rangle
\langle h_{1},h_2\rangle
+
\langle h_2,h_{1}^{\{ l_{ 1}-2\}}\rangle
h_2
<
$$
$$
<
\langle h_2,h_{1}^{\{ l_{ 1}-1\}}\rangle
\langle
h_{1}+1,h_2-1
\rangle
+
\langle h_2,h_{1}^{\{ l_{ 1}-2\}}\rangle
(h_2-1)
=
\langle h_2,h_{1}^{\{ l_{ 1}-1\}},
h_{1}+1,h_2-1
\rangle
.
$$
However, according
to the \ref{l11} Lemma, the inequality (\ref{vi9oh1})
contradicts the maximality of the continuant on the sequence $C=D_f$; therefore, $h_2=h_1+1$.

The rest of the proof follows from Lemma \ref{l18} and from Theorem \ref{t2}.
Theorem is proved.

\subsection{Proof  of Theorem \ref{t6}. }
Let $E$ --- be the finite sequence that gives the minimum value of the continuant on the set
$U(S, n,t)$, and let
$E\in W_f$.
Then by the Lemma \ref{l15} we get that $f\leq3$.
The minimal continuant of a finite sequence with at most three distinct elements ($1$, $n$, and $z$) can be found by the \ref{t3} theorem.
To apply the latter, imagine a sequence of$ t $ empty cells, which
we will begin to fill with the numbers $1$ and $n$ (alternating them), starting with the extreme empty cells, gradually moving to the central ones.
This is the reason for the formation of the central section $N_z (x)$ from (\ref{5}): if the elements ``$1$"\,  end first, then $N_z (x)$ will be as in the upper line from (\ref{5}), otherwise --- as in the lower one.
In this case, it is convenient to denote by $x$ or $-x$, respectively, the quantities ``$n$ " or "$1$"\,in $N_z (x)$ (except "$z$"), depending on which of the elements in $N_z(x)$ are present.
The rest of the proof follows from
the Lemma \ref{l19}.
The theorem is proved.

\subsection{Proof  of Theorem \ref{t7}. }
This minimum is also reduced to minimizing
the continuants of three-element permutations
by the \ref{t3} theorem, which results in a minimum in the form (\ref{12}).
But now the $t$ parameter is not fixed.
As in the proof of the previous theorem,
by $x$ we denote the number of ``$n$" \,in the final sequence, which will be as in the bottom line of (\ref{5}). 

Due to the inequality (\ref{9}), the condition $x\leq n$\,must be satisfied for the minimal continuant in the form (\ref{12}).
But since $x\geq 0$ and the continuant 
is taken from the sequence of the set $U_n (S)$, the number $S-z$
must be represented as the sum of several terms equal to 
 $n$ or $n+1$.
This is provided by selecting the set $P (S)$ in (\ref{7}), (\ref{8}).
Remaining relationship details 
(\ref{5})-(\ref{8}) and 
(\ref{19}), (\ref{20}) follow from the Lemma \ref{l20}. 
The theorem is proved.

\begin{Zam}\label{z4}
The extremuma from Theorems \ref{t1}-\ref{t6} are attained on the sequences which are unique up to
symmetry (\ref{31}) and extracting a unit  (\ref{32}) procedures.
It can be easily seen from the proofs.
\end{Zam}

\end{document}